\documentclass[12pt,reqno]{amsart}
\usepackage{amsfonts}
\usepackage{amsthm}
\usepackage{amsmath,amssymb}
\usepackage{amscd}
\usepackage[latin2]{inputenc}
\usepackage{t1enc}
\usepackage[mathscr]{eucal}
\usepackage{indentfirst}
\usepackage{graphicx}
\usepackage{graphics}
\usepackage{pict2e}
\usepackage{epic}
\usepackage[margin=2.9cm]{geometry}
\usepackage{epstopdf} 
\usepackage{bm}
\usepackage[all]{xy}
\usepackage{mathrsfs}
\usepackage{lmodern}
\newtheoremstyle{break}%
{}{}%
{\itshape}{}%
{\bfseries}{}
{\newline}{}

\newtheoremstyle{dbreak}%
{}{}%
{\upshape}{}%
{\bfseries}{}
{\newline}{}

\theoremstyle{break}
\newtheorem{Thm}{Theorem}[section]
\newtheorem{Lem}[Thm]{Lemma}
\newtheorem{Prop}[Thm]{Proposition}
\newtheorem{Cor}[Thm]{Corollary}

\newtheorem{theorem}{Theorem}

\theoremstyle{dbreak}

\newtheorem{Rem}[Thm]{Remark}
\newtheorem{Conj}[Thm]{Conjecture}

\newcommand{\prf}{\noindent\underline{$Proof.$}}

\numberwithin{equation}{section}

\newenvironment{nouppercase}{%
\renewcommand{\uppercasenonmath}[1]{}}{}

\makeatletter
\@namedef{subjclassname@2020}{%
  \textup{2020} Mathematics Subject Classification}
\makeatother

\begin{document}

\title[Anabelian properties of infinite algebraic extensions of finite fields]{Anabelian properties of infinite algebraic extensions of finite fields}
\author{Takahiro Murotani}

\subjclass[2020]{Primary 14H30; Secondary 11G20, 11S20, 14G15, 14G17, 14G20}

\keywords{anabelian geometry, \'{e}tale fundamental group, finite field, Grothendieck conjecture, higher local field, hyperbolic curve, infinite algebraic extension, Kummer-faithful field, tame fundamental group} 

\address[Takahiro Murotani]{Department of Mathematics, Tokyo Institute of Technology, Tokyo 152-8551, Japan}

\email{murotani.t.aa@m.titech.ac.jp}

\begin{abstract}
The Grothendieck conjecture for hyperbolic curves over finite fields was solved affirmatively by Tamagawa and Mochizuki.
On the other hand, (a ``weak version'' of) the Grothendieck conjecture for some hyperbolic curves over algebraic closures of finite fields is also known by Tamagawa and Sarashina.
So, it is natural to consider anabelian geometry over (infinite) algebraic extensions of finite fields.
In the present paper, we give certain generalizations of the above results of Tamagawa and Sarashina to hyperbolic curves over these fields.
Moreover, we give a necessary and sufficient condition for algebraic extensions of finite fields to be (torally) Kummer-faithful in terms of their absolute Galois groups.
\end{abstract}

\begin{nouppercase}
\maketitle
\end{nouppercase}

\tableofcontents

\section*{Introduction}

The Grothendieck conjecture for hyperbolic curves over finite fields was solved affirmatively by Tamagawa (the case where the curves are affine \cite[Theorems 0.5 and 0.6]{Tama1}) and Mochizuki (the general case \cite[Theorem 3.12]{cusp}).
(Some variants of these results are also known (cf. \cite[Theorem 1]{ST1} and \cite[Theorems C and D]{ST2}).)
On the other hand, (a ``weak version'' of) the Grothendieck conjecture for hyperbolic curves over algebraic closures of finite fields was partially solved by Tamagawa (the case where the genera of (the compactifications of) the curves are $0$ \cite[Theorem 0.2]{Tama2} and \cite[Theorem 0.2]{Tama3}) and Sarashina (the case where the curves are one-punctured elliptic curves and the characteristics are not equal to $2$ \cite[Theorem 1.2]{Sarashina}).

So, it is natural to consider anabelian geometry over infinite algebraic extensions of finite fields.
Some anabelian results such as the Grothendieck conjecture for hyperbolic curves over generalized sub-$p$-adic fields (cf. \cite[Theorem 4.12]{sur}), and criteria for the geometricity of homomorphisms between the absolute Galois groups of mixed-characteristic complete discrete valuation fields with perfect residue fields (cf. \cite[Theorems 0.1 and 0.2]{CHT}) seem to suggest the suitability of (infinite) algebraic extensions of finite fields for base fields of anabelian geometry.
(Note that these two results contain anabelian results over mixed-characteristic complete discrete valuation fields having these fields as residue fields.)
In the present paper, we consider the suitability of these fields for base fields of anabelian geometry from the following two points of view:
\begin{enumerate}
\item[(A)] Does the Grothendieck conjecture for hyperbolic curves over (infinite) algebraic extensions of finite fields hold?

\item[(B)] Can we describe the Kummer-faithfulness of (infinite) algebraic extensions of finite fields?
\end{enumerate}
We consider the problem (B) since Kummer-faithful fields are, to some extent, thought to be ``suitable'' fields by \cite[Theorem 3.4]{Kummer}.

\

For (A), we prove the following theorem:

\begin{theorem}[cf. Proposition \ref{topfin}, Corollaries \ref{type}, \ref{0}, \ref{wildtype} and \ref{wild0}]\label{(A)}
Let $k_i$ be an algebraic extension of a finite field, $\overline{k_i}$ an algebraic closure of $k_i$, $G_{k_i}$ the absolute Galois group $\mathrm{Gal}(\overline{k_i}/k_i)$, $X_i$ a hyperbolic curve over $k_i$, $X_i^{\mathrm{cpt}}$ the smooth compactification of $X_i$, $g_i$ the genus of $X_i^{\mathrm{cpt}}$, $S_i:=X_i^{\mathrm{cpt}}\setminus X_i$, $r_i$ the cardinality of $S_i(\overline{k_i})$, $\Pi_i$ either the \'{e}tale fundamental group $\pi_1(X_i)$ or the tame fundamental group $\pi_1^{\mathrm{tame}}(X_i)$ of $X_i$ and $\Delta_i$ the geometric subgroup of $\Pi_i$ for $i=1,\,2$ (for some choice of basepoint).
Suppose that there exists an isomorphism of profinite groups $\alpha:\Pi_1\stackrel{\sim}{\to}\Pi_2$.
Then the following assertions hold:

\begin{enumerate}
\item[(i)] $\Pi_1=\pi_1^{\mathrm{tame}}(X_1)$ if and only if $\Pi_2=\pi_1^{\mathrm{tame}}(X_2)$.

\item[(ii)] It holds that $\mathrm{char}\,k_1=\mathrm{char}\,k_2=:p$.
\end{enumerate}

In the following, suppose, moreover, that one of the following conditions holds:
\begin{enumerate}
\item[(I)] $X_i$ is proper for some (hence both of) $i=1,\,2$.

\item[(II)] The maximal pro-$p$ quotient $G_{k_i}^p$ of $G_{k_i}$ is trivial for some (hence both of) $i=1,\,2$.
\end{enumerate}
Then the following assertions hold:

\begin{enumerate}
\item[(iii)] It holds that $\alpha(\Delta_1)=\Delta_2$.
In particular, $\alpha$ induces an isomorphism $\overline{\alpha}:G_{k_1}\stackrel{\sim}{\to}G_{k_2}$ (in fact, $\overline{\alpha}$ is geometric (i.e., arises from a field isomorphism)).

\item[(iv)] $k_1$ and $k_2$ are isomorphic as fields.

\item[(v)] It holds that $(g_1,\,r_1)=(g_2,\,r_2)$.

\item[(vi)] Suppose that one of the following conditions holds:
\begin{enumerate}
\item[(a)] $g_1=g_2=0$.

\item[(b)] $(g_1,\,r_1)=(g_2,\,r_2)=(1,\,1)$ and $p\neq 2$.
\end{enumerate}
Then there exist finite extensions $k_1'$ and $k_2'$ of $k_1$ and $k_2$ (corresponding to each other via the isomorphism $\overline{\alpha}:G_{k_1}\stackrel{\sim}{\to}G_{k_2}$) such that $X_1\times_{\mathrm{Spec}\,k_1}\mathrm{Spec}\,k_1'$ and $X_2\times_{\mathrm{Spec}\,k_2}\mathrm{Spec}\,k_2'$ are isomorphic as schemes.
Moreover, if the condition (a) holds, we may take $k_1$ and $k_2$ as $k_1'$ and $k_2'$.

\item[(vii)] Suppose that the field $k_1\,(\simeq k_2)$ is finite.
Then $\alpha$ arises from a uniquely determined commutative diagram of schemes:
\[\xymatrix{\tilde{X_1} \ar[r]^{\sim} \ar[d] & \tilde{X_2} \ar[d] \\ X_1 \ar[r]^{\sim} & X_2}\]
in which the horizontal arrows are isomorphisms and the vertical arrows are the profinite \'{e}tale coverings corresponding to the groups $\Pi_1$ and $\Pi_2$, respectively.
\end{enumerate}
\end{theorem}
In fact, the conditions (I), (II) and the finiteness of $k_1\,(\simeq k_2)$ in (vii) of the above theorem are group-theoretic, i.e., we may determine from the isomorphism class of $\Pi_1\,(\simeq \Pi_2)$ whether these conditions hold or not (cf. Propositions \ref{topfin}, \ref{Gl}, \ref{ppart} and Remark \ref{tamemono}).

The key to the above theorem is to recover the projections from the \'{e}tale or tame fundamental groups to the absolute Galois groups of the base fields (cf. Proposition \ref{tame}, Lemmas \ref{prime} and \ref{ppart}).
This is possible at least in the case of tame fundamental groups and the case where the ``$p$-parts'' of the absolute Galois groups of the base fields are trivial (cf. the conditions (I) and (II) of Theorem \ref{(A)}).

\

On the other hand, for (B), we give a complete answer:

\begin{theorem}[cf. Proposition \ref{torally} and Theorem \ref{KF}]\label{(B)}
Let $k$ be an algebraic extension of a finite field of characteristic $p$ and $G_k$ the absolute Galois group of $k$.
\begin{enumerate}
\item[(i)] $k$ is torally Kummer-faithful if and only if the maximal pro-$l$ quotient $G_k^l$ of $G_k$ is isomorphic to $\mathbb{Z}_l$ for any $l\in\mathfrak{Primes}\setminus\{p\}$.

\item[(ii)] $k$ is Kummer-faithful if and only if $G_k$ is isomorphic to $\hat{\mathbb{Z}}$.
\end{enumerate}
\end{theorem}

By the above theorem, we can describe the (torally) Kummer-faithfulness of any algebraic extensions of finite fields.
It is a highly non-trivial problem to find the boundaries between Kummer-faithful fields and non-Kummer-faithful fields (i.e., to determine when infinite algebraic extensions of Kummer-faithful fields are again Kummer-faithful) (cf. \cite[\S 1]{Topics3} and \cite[\S 3]{OT}).
As far as the author knows, this is the first result describing such a boundary completely.

\

We shall review the contents of the present paper.
In Section 1, we investigate properties of the \'{e}tale and tame fundamental groups of hyperbolic curves over algebraic extensions of finite fields, and prove Theorem \ref{(A)}.
Moreover, we also give some ``pro-$\Sigma$ versions'' of anabelian results for hyperbolic curves over these fields, where $\Sigma$ is a set of prime numbers.
In Section 2, we prove Theorem \ref{(B)}, and consider properties of the absolute Galois groups of a certain type of higher local fields.

\

\section*{Acknowledgments}
The author would like to thank Professors Yuichiro Hoshi, Yoshiyasu Ozeki, Yuichiro Taguchi and Akio Tamagawa for constructive comments.
Especially, the author would like to express his gratitude to Professor Akio Tamagawa for informing him the arguments used to prove Corollary \ref{0} (and Corollary \ref{wild0}).
The author was supported by JSPS KAKENHI Grant Number 22J00022.

\

\

\setcounter{section}{-1}

\section{Notations and conventions}

\noindent
{\bf{Numbers and fields:}}

\

We shall write
\begin{itemize}
\item $\mathbb{Z}$ for the set of integers;

\item $\mathfrak{Primes}$ for the set of prime numbers.

\end{itemize}
For $a\in\mathbb{Z}$, we shall write $\mathbb{Z}_{\geq a}$ (resp. $\mathbb{Z}_{>a}$) for
\[\{b\in\mathbb{Z}\,|\,b\geq a\text{ (resp. $b>a$)}\}.\]

For $p\in\mathfrak{Primes}$ and $n\in\mathbb{Z}_{>0}$, we shall write
\begin{itemize}
\item $\mathbb{Z}_p$ for the $p$-adic completion of $\mathbb{Z}$;

\item $\mathbb{F}_{p^n}$ for the finite field of cardinality $p^n$.
\end{itemize}

\

\noindent
{\bf{Modules:}}

\

Let $M$ be a $\mathbb{Z}$-module.
We shall write $M_{\mathrm{div}}$ for the submodule which consists of divisible elements of $M$, i.e.,
\[M_{\mathrm{div}}=\bigcap_{n\in\mathbb{Z}_{>0}}nM.\]
For $n\in\mathbb{Z}_{>0}$, we shall write $M[n]$ for the kernel of multiplication by $n$, and set:
\[M[n^\infty]=\bigcup_{m\in\mathbb{Z}_{>0}}M[n^m].\]

\

\noindent
{\bf{Profinite groups:}}

\

Let $G$ be a profinite group, $p\in\mathfrak{Primes}$ and $\Sigma$ a set of prime numbers.
Then we shall write $G^\mathrm{ab}$ for the abelianization of $G$ (i.e., the quotient of $G$ by the closure of the commutator subgroup of $G$), $G^p$ (resp. $G^{p'}$, resp. $G^{\Sigma}$) for the maximal pro-$p$ quotient (resp. the maximal pro-prime-to-$p$ quotient, resp. the maximal pro-$\Sigma$ quotient) of $G$, $G^{\text{ab-tor}}$ for the subgroup consisting of elements of $G^{\text{ab}}$ of finite order, and $G^{\text{ab/tor}}$ for the quotient of $G^{\text{ab}}$ by the closure of $G^{\text{ab-tor}}$.

If $M$ is a topological $G$-module, we shall write $M_G$ for the maximal quotient of $M$ on which $G$ acts trivially.

We denote the cohomological $p$-dimension of $G$ by $\mathrm{cd}_pG$.

\

\

\section{Fundamental groups of hyperbolic curves over algebraic extensions of finite fields}

Let $p$ be a prime number and $k$ a(n) (not necessarily finite) algebraic extension of $\mathbb{F}_p$.
Moreover, we shall write:

\begin{itemize}
\item $\overline{k}$ for an algebraic closure of $k$;

\item $G_k$ for the Galois group $\mathrm{Gal}(\overline{k}/k)$;

\item $\mu_{l^\infty}\subset\overline{k}^\times$ for the group of $l$-power roots of unity ($l\in\mathfrak{Primes}\setminus\{p\}$);

\item $X$ for a hyperbolic curve over $k$;

\item $X^{\mathrm{cpt}}$ for the smooth compactification of $X$;

\item $g$ for the genus of $X^{\mathrm{cpt}}$;

\item $S:=X^{\mathrm{cpt}}\setminus X$;

\item $r$ for the cardinality of $S(\overline{k})$;

\item $\Pi_X$ for either the \'{e}tale fundamental group $\pi_1(X)$ of $X$ or the tame fundamental group $\pi_1^{\mathrm{tame}}(X)$ of $X$ (for some choice of basepoint);

\item $\Delta_X\subset\Pi_X$ for the quotient of the \'{e}tale fundamental group $\pi_1(X\times_{\mathrm{Spec}\,k}\mathrm{Spec}\,\overline{k})$ of $X_i\times_{\mathrm{Spec}\,k}\mathrm{Spec}\,\overline{k}$ determined by $\Pi_X$;

\item $J_X$ for the Jacobian variety of $X$;

\item $\gamma$ for the $p$-rank of $J_X$.
\end{itemize}

For a subset $\Sigma\subset\mathfrak{Primes}$, set $\Pi_X^{(\Sigma)}:=\Pi_X/\mathrm{Ker}(\Delta_X\twoheadrightarrow\Delta_X^\Sigma)$.

Moreover, set:
\[\varepsilon:=\begin{cases}1, \,&(r\neq 0), \\ 0, \,&(r=0).\end{cases}\]

Note that we have the following exact sequence:

\[\xymatrix{1 \ar[r] & \Delta_X^\Sigma \ar[r] & \Pi_X^{(\Sigma)} \ar[r]^{\mathrm{pr}_k} & G_k \ar[r] & 1.}\]

\

\begin{Prop}
Let $\Sigma\subset\mathfrak{Primes}$ be a non-empty subset.
Suppose that one of the following conditions holds:
\begin{enumerate}
\item[(i)] $\Sigma\neq\{p\}$.

\item[(ii)] $r\neq 0$.

\item[(iii)] $\gamma\neq 1$.
\end{enumerate}
Then $\Delta^\Sigma_X$ is strongly sn-internally indecomposable (in particular, slim) and sn-elastic.
(For the definition of sn-internally indecomposable (resp. slim, resp. sn-elastic) profinite groups, see, e.g., \cite[Definition 1.8]{MST} (resp. \cite[p.6]{MST}, resp. \cite[Definition 2.1]{MST}).)
Moreover, $\Pi_X^{(\Sigma)}$ is strongly sn-internally indecomposable (in particular, slim).
\end{Prop}

\prf

The assertion for $\Delta_X^\Sigma$ is \cite[Theorem 3.19]{MST}.
The assertion for $\Pi_X^{(\Sigma)}$ follows from this and \cite[Proposition 1.16]{MST}.
(For the slimness, see \cite[Remark 1.1.1]{MT2}.)
\qed

\

First, let us consider the case where $\Sigma=\mathfrak{Primes}$.
Similarly to \cite[Proposition 3.3 (i)]{Tama1}, the following holds:

\begin{Prop}\label{topfin}
The following conditions are all equivalent:

\begin{enumerate}
\item[(a)] Either $r=0$ or $\Pi_X=\pi_1^{\mathrm{tame}}(X)$.

\item[(b)] $\Pi_X$ is topologically finitely generated.

\item[(b$'$)] $\Pi_X^{\mathrm{ab}}$ is topologically finitely generated.

\item[(c)] $\Pi_X^p$ is topologically finitely generated.

\item[(c$'$)] $(\Pi_X^p)^{\mathrm{ab}}$ is topologically finitely generated.
\end{enumerate}
\end{Prop}

\

We may determine group-theoretically (from $\Pi_X$) whether $G_k^l$ is trivial or not for any $l\in\mathfrak{Primes}$ if $\Pi_X$ is topologically finitely generated (and for $l\in\mathfrak{Primes}\setminus\{p\}$ even if $\Pi_X$ is not topologically finitely generated):

\begin{Prop}\label{Gl}
Suppose that $\Pi_X$ is topologically finitely generated.
Then, for any prime number $l$, the following conditions are equivalent:

\begin{enumerate}
\item[(i)] $G_k^l\simeq\mathbb{Z}_l$.

\item[(ii)] For any open subgroup $H\subset\Pi_X$, $(H^l)^{\mathrm{ab\text{/}tor}}\simeq\mathbb{Z}_l$.
\end{enumerate}

Moreover, the above two conditions are equivalent for any prime number $l\neq p$ even in the case where $\Pi_X$ is not topologically finitely generated.
\end{Prop}

\prf

Set $\Pi_X^{(l)}:=\Pi_X/\mathrm{Ker}(\Delta_X\twoheadrightarrow\Delta_X^l)$.
Note that we have the following exact sequence:
\begin{align}\label{abq}
\xymatrix{1 \ar[r] & (\Delta_X^l)^{\mathrm{ab}}_{G_k} \ar[r] & (\Pi_X^{(l)})^{\mathrm{ab}} \ar[r] & G_k \ar[r] & 1.}
\end{align}

First, we shall prove that (i) implies (ii).
By replacing $X$ by a finite \'{e}tale covering, we assume that $H=\Pi_X$.
If $l\neq p$, we have the following exact sequence of $G_k$-modules (even in the case where $\Pi_X$ is not topologically finitely generated):
\begin{align}\label{Tate}
\xymatrix{0 \ar[r] & \mathbb{Z}_l(1)^{\oplus\varepsilon} \ar[r] & \mathbb{Z}[S(\overline{k})]\otimes_{\mathbb{Z}}\mathbb{Z}_l(1) \ar[r] & (\Delta_X^l)^{\mathrm{ab}} \ar[r] & T(J_X)^l \ar[r] & 0,}
\end{align}
where $T(J_X)$ is the Tate module of $J_X$.
Since we have assumed that (i) holds, we have $\mu_{l^\infty}\not\subset k^\times$.
Moreover, $T(J_X)^l_{G_k}$ is also finite.
Indeed, let $k_0$ be a finite subfield of $k$ over which $J_X$ is defined.
Since $G_k^l\simeq\mathbb{Z}_l$, $k$ is potential prime-to-$l$ extension of $k_0$ (for the definition of potential prime-to-$l$ extensions, see \cite[Definition 2.1]{Ozeki}).
Therefore, by \cite[Proposition 2.4]{Ozeki}, $J_X(k)[l^\infty]$ (hence also $T(J_X)^l_{G_k}$) is finite.
So, $(\Delta_X^l)^{\mathrm{ab}}_{G_k}$ is finite, and hence
\[(\Pi_X^l)^{\text{ab/tor}}=((\Pi_X^{(l)})^{\mathrm{ab}})^l)^{\text{ab/tor}}\simeq G_k^l\simeq\mathbb{Z}_l,\]
as desired.
On the other hand, if $l=p$ (and $\Pi_X$ is topologically finitely generated), we have $(\Delta_X^p)^{\mathrm{ab}}\simeq T(J_X)^p$.
Similarly as above, $J_X(k)[p^\infty]$ is finite and hence $T(J_X)^p_{G_k}$ is also finite.
Therefore, (ii) holds also in this case.

Next, suppose that (i) does not hold (i.e., $G_k^l$ is trivial).
Since $G_k$ is a projective profinite group, the exact sequence (\ref{abq}) splits.
Therefore, we have $(\Pi_X^l)^{\mathrm{ab}}\,(=((\Pi_X^{(l)})^{\mathrm{ab}})^l)\simeq (\Delta_X^l)^{\mathrm{ab}}_{G_k}$.
If $l\neq p$, there exists a finite extension $k'$ of $k$ such that $\mu_{l^\infty}\subset (k')^{\times}$ and that $J_X(k')[l^\infty]=J_X(\overline{k})[l^\infty]$ (cf. \cite[Proposition 2.2]{Ozeki}).
Then, for $H=\mathrm{pr}_k^{-1}(G_{k'})$, $(H^l)^{\mathrm{ab}}$ is isomorphic to free $\mathbb{Z}_l$-module of rank at least $2$.
If $l=p$, there also exists an open normal subgroup $H$ of $\Pi_X$ such that $(H^p)^{\mathrm{ab}}$ is isomorphic to free $\mathbb{Z}_p$-module of rank at least $2$ (cf. \cite[Lemma 1.9]{Tama1} and \cite[Proposition 2.2]{Ozeki}).
\qed

\

If $\Pi_X$ is topologically finitely generated, the subgroup $\Delta_X\subset\Pi_X$ is group-theoretically recovered:

\begin{Prop}\label{tame}
Suppose that $\Pi_X$ is topologically finitely generated.
Let $P$ be the set of all prime numbers $l$ satisfying one (and hence both) of two conditions in Proposition \ref{Gl}.
Then the following holds:
\[\Delta_X=\bigcap_{l\in P}\mathrm{Ker}(\Pi_X\twoheadrightarrow(\Pi_X^l)^{\mathrm{ab\text{/}tor}}).\]
\end{Prop}

\prf

Immediate from Proposition \ref{Gl} and its proof.
\qed

\begin{Rem}
In the case where $k$ is finite, Tamagawa recovered $\Delta_X$ from $\Pi_X$ group-theoretically (cf. \cite[Proposition 3.3 (ii)]{Tama1}).
However, the proof uses Weil conjecture, which we cannot apply in this case.
\end{Rem}

\begin{Rem}\label{tamemono}
Here, we give mono-anabelian reconstruction algorithms of various invariants of $X$ from $\Pi_X=\pi_1^{\mathrm{tame}}(X)$ explicitly (mainly according to \cite{Tama3}).

First, set:
\begin{align*}
P(\Pi_X)&:=\{l\in\mathfrak{Primes}\,|\,(H^l)^{\text{ab/tor}}\simeq\mathbb{Z}_l\text{ for any open subgroup $H\subset\Pi_X$}\}, \\
\Delta(\Pi_X)&:=\bigcap_{l\in P(\Pi_X)}\mathrm{Ker}(\Pi_X\twoheadrightarrow (\Pi_X^l)^{\text{ab/tor}}).
\end{align*}
Then, by Proposition \ref{tame}, it holds that $\Delta_X=\Delta(\Pi_X)$.

Next, let $p(\Pi_X)$ be the unique prime number $l$ satisfying the following condition:
\begin{quote}
For any prime number $l^\ast\neq l$, it holds that
\[\mathrm{rank}_{\mathbb{Z}_l}(\Delta(\Pi_X)^l)^{\mathrm{ab}}<\mathrm{rank}_{\mathbb{Z}_{l^\ast}}(\Delta(\Pi_X)^{l^\ast})^{\mathrm{ab}}.\]
\end{quote}
Moreover, set:
\[\gamma(\Pi_X):=\mathrm{rank}_{\mathbb{Z}_{p(\Pi_X)}}(\Delta(\Pi_X)^{p(\Pi_X)})^{\mathrm{ab}}.\]
Clearly, we have $p=p(\Pi_X),\,\gamma=\gamma(\Pi_X)$.

For a profinite group $G$ and $m\in\mathbb{Z}_{>0}$, we shall write $G(m)$ for the kernel of $G^{\mathrm{ab}}/m\cdot G^{\mathrm{ab}}$.
Fix a prime number $l\neq p(\Pi_X)$ and $m\in\mathbb{Z}_{>0}$ such that $m$ is prime to $p(\Pi_X)$ and $m\neq 1,\,3$, and set:
\begin{align*}
g'(\Pi_X)&:=\lim_{f\to\infty}\dfrac{\dim_{\mathbb{F}_{p(\Pi_X)}}((\Delta(\Pi_X)(p(\Pi_X)^f-1))^{\mathrm{ab}}\otimes\mathbb{F}_{p(\Pi_X)})}{\sharp(\Delta(\Pi_X)^{\mathrm{ab}}\otimes(\mathbb{Z}/(p(\Pi_X)^f-1)\mathbb{Z}))}=:g'(\Delta(\Pi_X)), \\
b(\Pi_X)&:=\mathrm{rank}_{\mathbb{Z}_l}\,H^1(\Delta(\Pi_X),\,\mathbb{Z}_l)=:b(\Delta(\Pi_X)), \\
r'(\Pi_X)&:=b(\Pi_X)-2g'(\Pi_X)+1=b(\Delta(\Pi_X))-2g'(\Delta(\Pi_X))+1=:r'(\Delta(\Pi_X)), \\
r(\Pi_X)&:=\begin{cases}r'(\Pi_X),\,&(r'(\Pi_X)\neq 3), \\ 0,\,&(r'(\Pi_X)=3,\,r'(\Delta(\Pi_X)(m))=3), \\ 1,\,&(r'(\Pi_X)=3,\,r'(\Delta(\Pi_X)(m))=m^{b(\Pi_X)}), \\ 3,\,&(r'(\Pi_X)=3,\,r'(\Delta(\Pi_X)(m))=3m^{b(\Pi_X)-1}),\end{cases} \\
g(\Pi_X)&:=\begin{cases}g'(\Pi_X)+1,\,&(r(\Pi_X)\leq 1), \\ g'(\Pi_X),\,&(r(\Pi_X)>1).\end{cases}
\end{align*}
(Note that the above invariants are well-defined and independent of the choices of $l$ and $m$ by arguments in \cite[\S 4]{Tama3}.)
Then, by arguments in \cite[\S 4]{Tama3}, it holds that $g=g(\Pi_X),\,r=r(\Pi_X)$.

Finally, set:
\begin{align*}
\mathcal{H}(\Pi_X)&:=\{H\subset\Pi_X\,|\,H\text{ is an open subgroup of $\Pi_X$},\,2g(H)-2=(\Pi_X:H)(2g(\Pi_X)-2)\}, \\
\Pi^\ast(\Pi_X)&:=\Pi_X/\bigcap_{H\in\mathcal{H}(\Pi_X)}H.
\end{align*}
By the Hurwitz formula, it holds that $(\Pi_X\twoheadrightarrow)\pi_1(X^\mathrm{cpt})\simeq\Pi^\ast(\Pi_X)$.
Let $U$ be an open subgroup of $\Pi_X$ such that $g(U)\geq 2$ (cf. \cite[Lemma 1.10]{Tama1}).
Then we have an open injection $U/\Delta(U)\hookrightarrow \Pi_X/\Delta(\Pi_X)$.
Set:
\begin{align*}
N(U)&:=(\Pi_X/\Delta(\Pi_X):U/\Delta(U)), \\
\Lambda(U)&:=\mathrm{Hom}_{\hat{\mathbb{Z}}^{p(\Pi_X)'}}(H^2(\Delta(\Pi^\ast(U)),\,\hat{\mathbb{Z}}^{p(\Pi_X)'}),\,\hat{\mathbb{Z}}^{p(\Pi_X)'}).
\end{align*}
Note that $\Lambda(U)$ is isomorphic to $\hat{\mathbb{Z}}^{p(\Pi_X)'}(1)$ as a $G_{k_U}\,(\simeq U/\Delta(U))$-module (cf., e.g., \cite[\S 1]{cusp}), where $k_U$ is the finite extension of $k$ determined by the open injection $U/\Delta(U)\hookrightarrow \Pi_X/\Delta(\Pi_X)\simeq G_k$.
Let $\overline{\mathbb{F}_{p(\Pi_X)}}$ be an algebraic closure of $\mathbb{F}_{p(\Pi_X)}$, and $\zeta_n\in\overline{\mathbb{F}_{p(\Pi_X)}}$ a primitive $n$-th root of unity for $n\in\mathbb{Z}_{>0}$ prime to $p(\Pi_X)$.
Set:
\[k'(U):=\mathbb{F}_{p(\Pi_X)}(\zeta_n\,|\,n\in\mathbb{Z}_{>0}\text{ is prime to $p(\Pi_X)$, and }U/\Delta(U)\text{ acts trivially on }\Lambda(U)/n\Lambda(U)),\]
and let $k(\Pi_X)$ be the unique intermediate field of $k'(U)/\mathbb{F}_p$ such that $[k'(U):k(\Pi_X)]=N(U)$.
Then $k(\Pi_X)$ is independent of the choice of $U$, and we have $k(\Pi_X)\simeq k$ as fields.
\end{Rem}

\begin{Cor}\label{type}
Let $k_i$ be an algebraic extension of a finite field, $\overline{k_i}$ an algebraic closure of $k_i$, $G_{k_i}$ the absolute Galois group $\mathrm{Gal}(\overline{k_i}/k_i)$, $X_i$ a hyperbolic curve over $k_i$, $X_i^{\mathrm{cpt}}$ the smooth compactification of $X_i$, $g_i$ the genus of $X_i^{\mathrm{cpt}}$, $S_i:=X_i^{\mathrm{cpt}}\setminus X_i$, $r_i$ the cardinality of $S_i(\overline{k_i})$, $\Pi_i$ the tame fundamental group $\pi_1^{\mathrm{tame}}(X_i)$ of $X_i$ and $\Delta_i$ the geometric subgroup of $\Pi_i$ for $i=1,\,2$ (for some choice of basepoint).
Suppose that there exists an isomorphism of profinite groups $\alpha:\Pi_1\stackrel{\sim}{\to}\Pi_2$.
Then the following assertions hold:

\begin{enumerate}
\item[(i)] It holds that $\alpha(\Delta_1)=\Delta_2$.
In particular, $\alpha$ induces an isomorphism $\overline{\alpha}:G_{k_1}\stackrel{\sim}{\to}G_{k_2}$ (in fact, $\overline{\alpha}$ is geometric (i.e., arises from a field isomorphism)).

\item[(ii)] $k_1$ and $k_2$ are isomorphic as fields.
In particular, it holds that $\mathrm{char}\,k_1=\mathrm{char}\,k_2=:p$.

\item[(iii)] It holds that $(g_1,\,r_1)=(g_2,\,r_2)$.

\item[(iv)] Suppose that one of the following conditions holds:
\begin{enumerate}
\item[(a)] $g_1=g_2=0$.

\item[(b)] $(g_1,\,r_1)=(g_2,\,r_2)=(1,\,1)$ and $p\neq 2$.
\end{enumerate}
Then there exist finite extensions $k_1'$ and $k_2'$ of $k_1$ and $k_2$ (corresponding to each other via the isomorphism $\overline{\alpha}:G_{k_1}\stackrel{\sim}{\to}G_{k_2}$) such that $X_1\times_{\mathrm{Spec}\,k_1}\mathrm{Spec}\,k_1'$ and $X_2\times_{\mathrm{Spec}\,k_2}\mathrm{Spec}\,k_2'$ are isomorphic as schemes.

\item[(v)] Suppose that the field $k(\Pi_1)\,(\simeq k(\Pi_2))$ constructed as in Remark \ref{tamemono} is finite.
Then $\alpha$ arises from a uniquely determined commutative diagram of schemes:
\[\xymatrix{\tilde{X_1} \ar[r]^{\sim} \ar[d] & \tilde{X_2} \ar[d] \\ X_1 \ar[r]^{\sim} & X_2}\]
in which the horizontal arrows are isomorphisms and the vertical arrows are the profinite \'{e}tale coverings corresponding to the groups $\Pi_1$ and $\Pi_2$, respectively.
\end{enumerate}
\end{Cor}

\prf

(i) is immediate from Proposition \ref{tame}.
(For the geometricity of $\overline{\alpha}$, see the claim in the proof of Corollary \ref{0}.)
(ii) and (iii) are immediate from Remark \ref{tamemono}.
(iv) is immediate from Proposition \ref{tame}, \cite[Theorems 0.1 and 0.2]{Tama3}, \cite[Theorem 1.2]{Sarashina} and their proofs.
(Though \cite[Theorem 1.2]{Sarashina} treats only the case of \'{e}tale fundamental groups, the same proof works also in the case of tame fundamental groups in light of arguments in \cite{Tama3} (especially, \cite[Theorems 4.1, 5.2 and Proposition 5.3]{Tama3}).)
(v) is immediate from \cite[Theorem 0.5]{Tama1} and \cite[Theorem 3.12]{cusp}.
\qed

\

In fact, in the genus $0$ case, we need not to take finite extensions of the base fields in the situation of (iv) of the above corollary:

\begin{Cor}\label{0}
In the situation of Corollary \ref{type}, suppose that $g_1=g_2=0$.
Then $X_1$ and $X_2$ are isomorphic as schemes.
\end{Cor}

\prf

First, we shall prove the following claim:

\noindent
\underline{Claim}

Let $K_i$ be an algebraic extension of a finite field of characteristic $p$, $\overline{K_i}$ an algebraic closure of $K_i$, $G_{K_i}$ the absolute Galois group $\mathrm{Gal}(\overline{K_i}/K_i)$ and $\chi_i:G_{K_i}\to(\hat{\mathbb{Z}}^{p'})^\times$ the cyclotomic character for $i=1,\,2$.
Suppose that there exists an isomorphism of profinite groups $\beta:G_{K_1}\stackrel{\sim}{\to}G_{K_2}$.
Then the following conditions are equivalent:
\begin{enumerate}
\item[(i)] $\beta$ is geometric (i.e., arises from a field isomorphism $k_1\simeq k_2$).

\item[(ii)] It holds that $\phi_1=\phi_2\circ\beta$, where $\phi_i:G_{K_i}\to G_{\mathbb{F}_p}$ is induced by the field embedding $\mathbb{F}_p\hookrightarrow K_i$ ($i=1,\,2$) and $G_{\mathbb{F}_p}$ the absolute Galois group of $\mathbb{F}_p$.

\item[(iii)] It holds that $\chi_1^n=\chi_2^n\circ\beta$ for some $n\in\mathbb{Z}_{>0}$.
\end{enumerate}
The implications (i)$\Longrightarrow$(ii)$\Longrightarrow$(iii) are immediate.
Suppose that (iii) holds.
Let $\chi:G_{\mathbb{F}_p}\to(\hat{\mathbb{Z}}^{p'})^\times$ be the cyclotomic character.
Then it holds that $\chi^n\circ\phi_1=\chi^n\circ\phi_2\circ\beta$.
However, since $\chi^n$ is injective, we have $\phi_1=\phi_2\circ\beta$, hence (ii) holds.
Next, suppose that (ii) holds.
Then $K_1$ is isomorphic to $K_2$ as a(n) (abstract) field.
So we may assume that $K_1=K_2$.
However, in this case, $\beta$ must be the identity, hence geometric, as desired.
This completes the proof of the claim.

By Proposition \ref{tame}, $\alpha$ determines an isomorphism $\Delta_1\stackrel{\sim}{\to}\Delta_2$.
Set $r:=r_1=r_2$ (cf. Corollary \ref{type} (iii)).
By the arguments in \cite[\S 4, 5]{Tama3} (see also \cite[\S 3]{Tama2}), there exist a field isomorphism $\iota:\overline{k_1}\stackrel{\sim}{\to}\overline{k_2}$ and $\lambda_1,\,\cdots,\,\lambda_{r-3}\in\overline{k_1},\,\mu_1,\,\cdots,\,\mu_{r-3}\in\overline{k_2}$ such that
\begin{align*}
(X_1)_{\overline{k_1}}:=X_1\times_{\mathrm{Spec}\,k_1}\mathrm{Spec}\,\overline{k_1}&\simeq\mathbb{P}^1_{\overline{k_1}}\setminus\{0,\,1,\,\infty,\,\lambda_1,\,\cdots,\,\lambda_{r-3}\}, \\
(X_2)_{\overline{k_2}}:=X_2\times_{\mathrm{Spec}\,k_2}\mathrm{Spec}\,\overline{k_2}&\simeq\mathbb{P}^1_{\overline{k_1}}\setminus\{0,\,1,\,\infty,\,\mu_1,\,\cdots,\,\mu_{r-3}\},
\end{align*}
and that $\iota(\lambda_j)=\mu_j\,(1\leq j\leq r-3)$.
(Note that, since $k_1\simeq k_2$ (cf. Corollary \ref{type} (iv)), $\iota$ must map $k_1$ into $k_2$.)

On the other hand, by a similar argument in the proof of \cite[Proposition 3.4 (i)]{Tama1}, it holds that $\chi_1^n=\chi_2^n\circ\overline{\alpha}$ for some $n\in\mathbb{Z}_{>0}$ (e.g., $n=2(r-1)$ (note that, in this case, $r\geq 3$)).
By the above claim, $\overline{\alpha}$ is geometric and hence induced by a field isomorphism $\iota':\overline{k_1}\stackrel{\sim}{\to}\overline{k_2}$ (which maps $k_1$ into $k_2$).
Then there exists an element $\tau\in\mathrm{Aut}_{\mathbb{F}_p}(\overline{k_2})=\mathrm{Gal}(\overline{k_2}/\mathbb{F}_p)$ such that $\iota=\tau\circ\iota'$.
By definition, $\iota'$ is Galois equivariant with respect to $\overline{\alpha}$ (i.e., for any $\sigma\in G_{k_1}$, it holds that $\iota'\circ\sigma=\overline{\alpha}(\sigma)\circ\iota'$).
Since $\mathrm{Gal}(\overline{k_2}/\mathbb{F}_p)\,(\simeq\hat{\mathbb{Z}})$ is abelian, $\iota$ is also Galois equivariant (with respect to $\overline{\alpha}$) in this sense.
Moreover, by \cite[\S 2]{Tama2}, there exists a bijection $f:(S_1)_{\overline{k_1}}\stackrel{\sim}{\to}(S_2)_{\overline{k_2}}$ which is compatible with the Galois actions and $\overline{\alpha}$.
Therefore, we have the following commutative diagram:

\[\xymatrix{(S_1)_{\overline{k_1}} \ar@{}[r]|{\subset} \ar[d]^\wr_f & (X_1^{\mathrm{cpt}})_{\overline{k_1}} \ar[r] \ar[d]^\wr & \mathrm{Spec}\,\overline{k_1} \ar[d]^\wr_{\mathrm{Spec}\,(\iota^{-1})}\\ (S_2)_{\overline{k_2}} \ar@{}[r]|{\subset}  & (X_2^{\mathrm{cpt}})_{\overline{k_2}} \ar[r] & \mathrm{Spec}\,\overline{k_2}}\]
Since $f$ and $\mathrm{Spec}\,(\iota^{-1})$ is compatible with Galois actions with respect to $\overline{\alpha}$, the isomorphism $(X_1^{\mathrm{cpt}})_{\overline{k_1}}\stackrel{\sim}{\to}(X_2^{\mathrm{cpt}})_{\overline{k_2}}$ is also compatible with Galois actions (note that an automorphism of a projective line is determined by the images of three distinct points).
By taking quotients by the Galois actions, we have
\[\xymatrix{S_1 \ar@{}[r]|{\subset} \ar[d]^\wr&  X_1^{\mathrm{cpt}} \ar[d]^\wr \\ S_2 \ar@{}[r]|{\subset} & X_2^{\mathrm{cpt}}}\]
and hence an isomorphism $X_1=X_1^{\mathrm{cpt}}\setminus S_1\simeq X_2^{\mathrm{cpt}}\setminus S_2=X_2$.
\qed

\

\

Next, let us consider the case where $\Pi_X$ is not topologically finitely generated.
Note that the characteristic $p$ of $k$ is group-theoretically recovered from $\Pi_X$ as the unique prime number $l$ such that $\Pi_X^l$ is not topologically finitely generated (cf. Proposition \ref{topfin}).

Also in this case, a similar result to Proposition \ref{tame} holds for $l\in\mathfrak{Primes}\setminus\{p\}$:

\begin{Lem}\label{prime}
Suppose that $\Pi_X$ is not topologically finitely generated.
Let $P$ be the set of all prime numbers $l\neq p$ satisfying one (and hence both) of two conditions in Proposition \ref{Gl}.
Then the following holds:
\[\mathrm{Ker}(\Pi_X\twoheadrightarrow G_k^{p'})=\bigcap_{l\in P}\mathrm{Ker}(\Pi_X\twoheadrightarrow(\Pi_X^l)^{\mathrm{ab\text{/}tor}}).\]
\end{Lem}

\prf

Immediate from Proposition \ref{Gl} and its proof.
\qed

\

Moreover, we may determine group-theoretically whether $G_k^p$ is trivial or not:

\begin{Lem}\label{ppart}
Suppose that $\Pi_X$ is not topologically finitely generated.
Let $k'$ be the maximal prime-to-$p$ extension of $k$ (in $\overline{k}$) and $\Pi':=\mathrm{pr}_k^{-1}(G_{k'})$ (note that $\Pi'=\mathrm{Ker}(\Pi_X\twoheadrightarrow G_k^{p'})$ (cf. Lemma \ref{prime})).
Then the following conditions are equivalent:

\begin{enumerate}
\item[(i)] $G_k^p\,(\simeq G_{k'}^p)$ is trivial (i.e., $k'=\overline{k}$).

\item[(ii)] For any open subgroup $H$ of $\Pi'$, any open subgroup $H'$ of $H$ and any prime number $l\neq p$, it holds that
\[\mathrm{rank}_{\mathbb{Z}_l}((H')^l)^{\mathrm{ab\text{/}tor}}>\mathrm{rank}_{\mathbb{Z}_l}(H^l)^{\mathrm{ab\text{/}tor}}.\]
\end{enumerate}
\end{Lem}

\prf

By replacing $X$ by $X\times_{\mathrm{Spec}\,k}\mathrm{Spec}\,k'$ if necessary, we assume that $\Pi_X=\Pi'$ and $k=k'$.

First, suppose that (i) holds.
Fix a prime number $l\neq p$ and let $Y$ (resp. $Y'$) be the finite \'{e}tale covering of $X$ corresponding to $H$ (resp. $H'$), $g_Y$ (resp. $g_{Y'}$) the genus of the smooth compactification $Y^{\mathrm{cpt}}$ (resp. $(Y')^{\mathrm{cpt}}$) of $Y$ (resp. $Y'$) and $r_Y$ (resp. $r_{Y'}$) the cardinality of $(Y^{\mathrm{cpt}}\setminus Y)(\overline{k})$ (resp. $((Y')^{\mathrm{cpt}}\setminus Y')(\overline{k})$).
By the Hurwitz formula, we have
\[\mathrm{rank}_{\mathbb{Z}_l}((H')^l)^{\mathrm{ab\text{/}tor}}=2g_{Y'}+r_{Y'}-1>2g_Y+r_Y-1=\mathrm{rank}_{\mathbb{Z}_l}(H^l)^{\mathrm{ab\text{/}tor}},\]
as desired.

Next, suppose that (i) does not hold, and we shall prove that (ii) does not hold.
Fix a prime number $l\neq p$.
By \cite[Proposition 2.2]{Ozeki}, there exists a finite extension $k_1$ of $k$ such that $J_X(k_1)[l^\infty]=J_X(\overline{k})[l^\infty]$ and that $\mu_{l^\infty}\subset k_1^\times$.
Let $k_2$ be the subextension of $\overline{k}/k_1$ of degree $p$ and set $H:=\mathrm{pr}_k^{-1}(G_{k_1}),\,H':=\mathrm{pr}_k^{-1}(G_{k_2})$.
Then, clearly we have
\[\mathrm{rank}_{\mathbb{Z}_l}((H')^l)^{\mathrm{ab\text{/}tor}}=\mathrm{rank}_{\mathbb{Z}_l}(H^l)^{\mathrm{ab\text{/}tor}}=2g+r-1,\]
as desired.
(Note that we have $(H^l)^{\mathrm{ab\text{/}tor}}\simeq(\Delta_X^l)^{\mathrm{ab}}_{G_{k_1}}=(\Delta_X^l)^{\mathrm{ab}}$ (cf. the proof of Proposition \ref{Gl}).)
\qed

\

Although the author at the time of writing does not know how to recover the subgroup $\mathrm{Ker}(\Pi_X\twoheadrightarrow G_k^p)$ group-theoretically in the case where $\Pi_X$ is not topologically finitely generated, we might conjecture:

\begin{Conj}
Suppose that $\Pi_X$ is not topologically finitely generated, and that $G_k\simeq\mathbb{Z}_p$ (cf. Lemmas \ref{prime} and \ref{ppart}).
Then $\Delta_X$ is characterized as the unique closed normal subgroup $N$ of $\Pi_X$ satisfying the following conditions:

\begin{enumerate}
\item[(i)] $\Pi_X/N\simeq\mathbb{Z}_p$.

\item[(ii)] For $m\in\mathbb{Z}_{>0}$, set $N_m:=\mathrm{Ker}(\Pi_X\twoheadrightarrow (\Pi_X/N)/p^m(\Pi_X/N))$.
Moreover, for $m\in\mathbb{Z}_{>0}$, an open subgroup $U$ of $\Pi_X$ and $l\in\mathfrak{Primes}\setminus\{p\}$, set:
\[R_m(U,\,l):=\mathrm{rank}_{\mathbb{Z}_l}\,((U\cap N_m)^l)^{\mathrm{ab\text{/}tor}}.\]
Then, for any open subgroup $U$ of $\Pi_X$, there exists a positive integer $R(U)$ (depending only on $U$) such that, for any $l\in\mathfrak{Prime}\setminus\{p\}$, it holds that
\[R(U)=\lim_{m\to\infty}R_m(U,\,l).\]
\end{enumerate}
\end{Conj}

\begin{Rem}\label{wildmono}
Here, we give mono-anabelian reconstruction algorithms of various invariants of $X$ from $\Pi_X=\pi_1(X)$ explicitly in the case where $\Pi_X$ is not topologically finitely generated (mainly according to \cite{Tama2}).
For the case where $\Pi_X$ is topologically finitely generated, see Remark \ref{tamemono}.

First, let $p(\Pi_X)$ be the unique prime number $l$ such that $\Pi_X^l$ is not topologically finitely generated.
Then we have $p=p(\Pi_X)$ (cf. the paragraph preceding Lemma \ref{prime}).
Set:
\begin{align*}
P(\Pi_X)&:=\{l\in\mathfrak{Primes}\setminus\{p(\Pi_X)\}\,|\,(H^l)^{\text{ab/tor}}\simeq\mathbb{Z}_l\text{ for any open subgroup $H\subset\Pi_X$}\}, \\
\Delta'(\Pi_X)&:=\bigcap_{l\in P(\Pi_X)}\mathrm{Ker}(\Pi_X\twoheadrightarrow (\Pi_X^l)^{\text{ab/tor}}).
\end{align*}
In the following, suppose that the condition (ii) of Lemma \ref{ppart} for $\Pi'=\Delta'(\Pi_X)$ and $p=p(\Pi_X)$ holds.
Then, by Lemmas \ref{prime} and \ref{ppart}, it holds that $\Delta_X=\Delta'(\Pi_X)$.

Set:
\begin{align*}
\chi(\Pi_X)&:=1-\mathrm{rank}_{\hat{\mathbb{Z}}^{p(\Pi_X)'}}((\Delta'(\Pi_X)^{p(\Pi_X)'})^{\mathrm{ab}})\,(=:\chi(\Delta'(\Pi_X))), \\
\mathcal{H}(\Pi_X)&:=\{H\subset\Delta'(\Pi_X)\,|\,H\text{ is an open subgroup of $\Delta'(\Pi_X)$},\,\chi(H)=(\Delta'(\Pi_X):H)\chi(\Pi_X)\} \\
(&=:\mathcal{H}(\Delta'(\Pi_X))), \\
\Delta^t(\Pi_X)&:=\Delta'(\Pi_X)/\bigcap_{H\in\mathcal{H}(\Pi_X)}H\,(=:\Delta^t(\Delta'(\Pi_X))).
\end{align*}
Then, by \cite[Corollary 1.5]{Tama2}, we have $(\Delta_X\twoheadrightarrow)\pi_1^{\mathrm{tame}}(X\times_{\mathrm{Spec}\,k}\mathrm{Spec}\,\overline{k})\simeq\Delta^t(\Pi_X)$.
Moreover, set:
\[\gamma(\Pi_X):=\mathrm{rank}_{\mathbb{Z}_{p(\Pi_X)}}(\Delta^t(\Pi_X)^{p(\Pi_X)})^{\mathrm{ab}}\,(=:\gamma(\Delta'(\Pi_X))).\]
By \cite[Corollaries 1.7 and 1.8]{Tama2}, we have $\gamma=\gamma(\Pi_X)$.
Set:
\begin{align*}
r(\Pi_X)&:=\dfrac{1}{p(\Pi_X)-1}\max_{H}(\gamma(H)-1-p(\gamma(\Pi_X)-1)), \\
g(\Pi_X)&:=\dfrac{1}{2}(2-\chi(\Pi_X)-r(\Pi_X)),
\end{align*}
where $H$ runs through the set of open normal subgroups of $\Delta'(\Pi_X)$ of index $p(\Pi_X)$.
Then, by \cite[Theorem 1.9]{Tama2}, we have $r=r(\Pi_X),\,g=g(\Pi_X)$.

Finally, by taking an open subgroup $U$ of $\Pi_X$ such that $g(U)\geq 2$, we may construct a field $k(\Pi_X)$ isomorphic to $k$ by a similar argument to Remark \ref{tamemono}.
\end{Rem}

\begin{Cor}\label{wildtype}
Let $k_i$ be an algebraic extension of a finite field, $\overline{k_i}$ an algebraic closure of $k_i$, $G_{k_i}$ the absolute Galois group $\mathrm{Gal}(\overline{k_i}/k_i)$, $X_i$ be a hyperbolic curve over $k_i$, $X_i^{\mathrm{cpt}}$ the smooth compactification of $X_i$, $g_i$ the genus of $X_i^{\mathrm{cpt}}$, $S_i:=X_i^{\mathrm{cpt}}\setminus X_i$, $r_i$ the cardinality of $S_i(\overline{k_i})$, $\Pi_i$ the \'{e}tale fundamental group $\pi_1(X_i)$ of $X_i$ and $\Delta_i$ the geometric subgroup of $\Pi_i$ for $i=1,\,2$ (for some choice of basepoint).
Suppose that there exists an isomorphism of profinite groups $\alpha: \Pi_1\stackrel{\sim}{\to}\Pi_2$.
Then the following assertion holds:

\begin{enumerate}
\item[(i)] It holds that $\mathrm{char}\,k_1=\mathrm{char}\,k_2=:p$.

\end{enumerate}

In the following, suppose, moreover, that one of the following conditions holds:
\begin{enumerate}
\item[(a)] $X_i$ is proper for some (hence both of) $i=1,\,2$ (cf. Proposition \ref{topfin}).

\item[(b)] $G_{k_i}^p$ is trivial for some (hence both of) $i=1,\,2$ (cf. Proposition \ref{Gl} and Lemma \ref{ppart}).
\end{enumerate}
Then the following assertions hold:

\begin{enumerate}
\item[(ii)] It holds that $\alpha(\Delta_1)=\Delta_2$.
In particular, $\alpha$ induces an isomorphism $\overline{\alpha}:G_{k_1}\stackrel{\sim}{\to}G_{k_2}$ (in fact, $\overline{\alpha}$ is geometric (i.e., arises from a field isomorphism)).

\item[(iii)] $k_1$ and $k_2$ are isomorphic as fields.

\item[(iv)] It holds that $(g_1,\,r_1)=(g_2,\,r_2)$.

\item[(v)] Suppose that one of the following conditions holds:
\begin{enumerate}
\item[(a)] $g_1=g_2=0$.

\item[(b)] $(g_1,\,r_1)=(g_2,\,r_2)=(1,\,1)$ and $p\neq 2$.
\end{enumerate}
Then there exist finite extensions $k_1'$ and $k_2'$ of $k_1$ and $k_2$ (corresponding to each other via the isomorphism $\overline{\alpha}:G_{k_1}\stackrel{\sim}{\to}G_{k_2}$) such that $X_1\times_{\mathrm{Spec}\,k_1}\mathrm{Spec}\,k_1'$ and $X_2\times_{\mathrm{Spec}\,k_2}\mathrm{Spec}\,k_2'$ are isomorphic as schemes.
\end{enumerate}

\end{Cor}

\prf

(i) is immediate from Remark \ref{tamemono} and the paragraph preceding Lemma \ref{prime} (see also Remark \ref{wildmono}).
(ii) is immediate from Proposition \ref{tame} and Lemma \ref{prime}.
(For the geometricity of $\overline{\alpha}$, see the claim in the proof of Corollary \ref{0}.)
(iii) and (iv) are immediate from Remarks \ref{tamemono} and \ref{wildmono}.
(v) is immediate from Lemma \ref{ppart}, \cite[Theorem 0.2]{Tama2}, \cite[Theorem 1.2]{Sarashina} and their proofs.
\qed

\

Similarly to Corollary \ref{0}, in the genus $0$ case, we obtain:

\begin{Cor}\label{wild0}
In the situation of Corollary \ref{wildtype}, suppose that the condition (b) (of Corollary \ref{wildtype}) is satisfied, and that $g_1=g_2=0$.
Then $X_1$ and $X_2$ are isomorphic as schemes.
\end{Cor}

\

\

Finally, let us consider the case where $\Sigma\subsetneq\mathfrak{Primes}$.
A similar result to Proposition \ref{Gl} holds, i.e., the following conditions are equivalent for any $l\in\mathfrak{Primes}\setminus\{p\}$ (and for $l=p$ if either $p\not\in\Sigma$ or $\Pi_X=\pi_1^{\mathrm{tame}}(X)$):
\begin{enumerate}
\item[(i)] $G_k^l\simeq\mathbb{Z}_l$.

\item[(ii)] For any open subgroup $H\subset\Pi_X^{(\Sigma)}$, $(H^l)^{\mathrm{ab\text{/}tor}}\simeq\mathbb{Z}_l$.
\end{enumerate}
Therefore, similarly to Proposition \ref{tame}, the subgroup $\Delta_X^\Sigma\subset\Pi_X^{(\Sigma)}$ is recovered group-theoretically if either $p\not\in\Sigma$ or $\Pi_X=\pi_1^{\mathrm{tame}}(X)$:
\[\Delta_X^\Sigma=\bigcap_{l\in P}\mathrm{Ker}(\Pi_X^{(\Sigma)}\twoheadrightarrow ((\Pi_X^{(\Sigma)})^l)^{\text{ab/tor}}),\]
where $P$ is the set of all prime numbers $l$ satisfying one (and hence both) of above two conditions.

In general, it seems to be impossible to recover $p$, $g$, $r$ or the isomorphism class of $k$ group-theoretically.
However, in the following, consider the case where the following conditions hold:
\begin{enumerate}
\item[(a)] Either $p\not\in\Sigma$ or $\Pi_X=\pi_1^{\mathrm{tame}}(X)$.

\item[(b)] The set $\mathfrak{Primes}\setminus\Sigma$ is finite.
\end{enumerate}
Note that we may determine group-theoretically from $\Pi_X^{(\Sigma)}$ whether the above conditions (a) and (b) hold or not.
Indeed, (a) holds if and only if $(\Pi_X^{(\Sigma)})^l$ is topologically finitely generated for any $l\in\mathfrak{Primes}$.
Moreover, under the assumption that the condition (a) holds (hence $\Delta_X^\Sigma\subset\Pi_X^{(\Sigma)}$ is already recovered group-theoretically), (b) holds if and only if the set $\{l\in\mathfrak{Primes}\,|\,(\Delta_X^\Sigma)^l\text{ is trivial}\}$ is finite.

Under these assumptions, the following conditions are equivalent:
\begin{enumerate}
\item[(I)] $k$ is finite.

\item[(II)] $\Pi_X^{(\Sigma)}/\Delta_X^\Sigma\,(\simeq G_k)\simeq\hat{\mathbb{Z}}$ and $(\Pi_X^{(\Sigma)})^{\text{ab-tor}}$ is finite.
\end{enumerate}
The implication (I)$\Longrightarrow$(II) is immediate (cf. the proof of Proposition \ref{Gl}).
Suppose that (I) does not hold, and we shall show that (II) does not hold.
We may assume that $\Pi_X^{(\Sigma)}/\Delta_X^\Sigma\simeq G_k\simeq\hat{\mathbb{Z}}$.
Since $k$ is infinite, $T(J_X)_{G_k}$ is infinite.
On the other hand, since we have $G_k\simeq\hat{\mathbb{Z}}$, as in the proof of Proposition \ref{Gl}, $J_X(k)[l^\infty]$ is finite for any $l\in\mathfrak{Primes}$.
This shows that $T(J_X)^\Sigma_{G_k}$ is infinite.
Therefore, $(\Pi_X^{(\Sigma)})^{\text{ab-tor}}=(\Delta_X^\Sigma)^{\text{ab-tor}}_{G_k}$ is also infinite, as desired.
So, under the assumptions (a) and (b), we may determine group-theoretically (from $\Pi_X^{(\Sigma)}$) whether $k$ is finite or not.

In the case where $k$ is finite (and the above conditions (a) and (b) hold), we may give mono-anabelian reconstruction algorithms of the Frobenius element of $G_k$, $p$, $g$, $r$, the isomorphism class of $k$ and the subgroup $\mathrm{Ker}(\Pi_X^{(\Sigma)}\twoheadrightarrow\Pi_{X^{\mathrm{cpt}}}^{(\Sigma)})\,(=\mathrm{Ker}(\Delta_X^\Sigma\twoheadrightarrow\Delta_{X^{\mathrm{cpt}}}^\Sigma))$ (from $\Pi_X^{(\Sigma)}$) (cf. \cite[Proposition 1.15]{ST1} and \cite[\S 3]{Tama1}).
Moreover, by applying \cite[Theorems C and D]{ST2}, we may recover the isomorphism class of $X$ in this case:

\begin{Thm}\label{refined}
Let $k_i$ be an algebraic extension of a finite field, $p_i$ the characteristic of $k_i$, $X_i$ a hyperbolic curve over $k_i$, $\Sigma_i\subset\mathfrak{Primes}$ a set of prime numbers and $\Pi_i$ the geometrically pro-$\Sigma_i$ tame fundamental group of $X_i$ (i.e., $\Pi_i=\pi_1^{\mathrm{tame}}(X_i)/\mathrm{Ker}(\Delta_{X_i}^{\mathrm{tame}}\twoheadrightarrow(\Delta_{X_i}^{\mathrm{tame}})^{\Sigma_i})$, where $\pi_1^{\mathrm{tame}}(X_i)$ is the tame fundamental group of $X_i$ and $\Delta_{X_i}^{\mathrm{tame}}$ the geometric subgroup of $\pi_1^{\mathrm{tame}}(X_i)$ (for some choice of basepoint)) for $i=1,\,2$.
Suppose that there exists an isomorphism of profinite groups $\alpha:\Pi_1\stackrel{\sim}{\to}\Pi_2$, and that the following conditions hold:

\begin{enumerate}
\item[(i)] $\mathfrak{Primes}\setminus\Sigma_i$ is finite for some (hence both of) $i=1,\,2$ (cf. the arguments preceding this theorem).

\item[(ii)] $G_{k_i}\simeq\hat{\mathbb{Z}}$ for some (hence both of) $i=1,\,2$ (cf. the arguments preceding this theorem).

\item[(iii)] $(\Pi_1)^{\mathrm{ab\text{-}tor}}\,(\simeq(\Pi_2)^{\mathrm{ab\text{-}tor}})$ is finite.

\end{enumerate}
Then $k_1$ and $k_2$ are finite fields, and isomorphic as fields.
Moreover, $\alpha$ arises from a uniquely determined commutative diagram of schemes:
\[\xymatrix{\tilde{X_1} \ar[r]^{\sim} \ar[d] & \tilde{X_2} \ar[d] \\ X_1 \ar[r]^{\sim} & X_2}\]
in which the horizontal arrows are isomorphisms and the vertical arrows are the profinite \'{e}tale coverings corresponding to the groups $\Pi_1$ and $\Pi_2$, respectively.
\end{Thm}

\prf

Immediate from \cite[Theorems C and D]{ST2} and the arguments preceding this theorem.
\qed

\begin{Rem}
If we already know the characteirstic $p$ of $k$ and that $k$ is finite, we may recover the isomorphism class of $X$ from $\Pi_X^{(\Sigma)}$ under the assumption that the $\Sigma'$-adic representation $\rho_{\Sigma'}: G_k\to\displaystyle\prod_{l\in\Sigma'}\mathrm{GL}(T(J_X)^l)$ is not injective (cf. \cite[Theorem C]{ST2}), where $\Sigma':=\mathfrak{Primes}\setminus\Sigma$.
However, in the case where $\Sigma'$ is not necessarily finite, the author at the time of writing does not know how to recover $p\,(=\mathrm{char}\,k)$ group-theoretically from $\Pi_X^{(\Sigma)}$.
(Note that \cite[Proposition 2.5]{ST2} is applicable only when we already know that $k$ is finite.)
\end{Rem}

\

\

\section{The Kummer-faithfulness of fields algebraic over finite fields}\label{CHT}

In the present section, we give a characterization of algebraic extensions of finite fields which are (torally) Kummer-faithful.
(For the definition of Kummer-faithful (resp. torally Kummer-faithful) fields, see, e.g., \cite[Definition 2.2]{OT} (resp. \cite[\S 3.2]{OT}).)

\

Let $k$ be a(n) (not necessarily finite) algebraic extension of $\mathbb{F}_p$.
Moreover, we shall write:

\begin{itemize}
\item $\overline{k}$ for an algebraic closure of $k$;

\item $G_k$ for the Galois group $\mathrm{Gal}(\overline{k}/k)$;

\item $\mu_{l^\infty}\subset\overline{k}^\times$ for the group of $l$-power roots of unity ($l\in\mathfrak{Primes}\setminus\{p\}$).
\end{itemize}

\

\begin{Prop}\label{torally}
$k$ is torally Kummer-faithful if and only if $G_k^l\simeq\mathbb{Z}_l$ for any $l\in\mathfrak{Primes}\setminus\{p\}$.
\end{Prop}

\prf

First, suppose that $k$ is torally Kummer-faithful.
Since $k$ is algebraic over $\mathbb{F}_p$, $G_k\simeq\displaystyle\prod_{l\in P}\mathbb{Z}_l$ for some (possibly empty) set $P$ of prime numbers.
Suppose that $P$ does not contain a prime number $l_0\in\mathfrak{Primes}\setminus\{p\}$.
By replacing  $k$ by a finite extension, we assume that $k$ contains a primitive $l_0$-th root of unity $\zeta_{l_0}$.
Then, by the assumption that $G_k^{l_0}$ is trivial, it holds that $\mu_{l_0^\infty}\subset k^\times$.
Therefore, $\zeta_{l_0}\in(k^\times)_{\mathrm{div}}$.
This contradicts the assumption that $k$ is torally Kummer-faithful.

Conversely, suppose that $k$ is not torally Kummer-faithful.
By replacing $k$ by a finite extension if necessary, we assume that $k^\times$ contains a non-trivial divisible element $a\in k^\times$.
Since $k$ is algebraic over $\mathbb{F}_p$, $a$ is a primitive $N$-th root of unity for some integer $N\,(>1)$ which is coprime to $p$.
Let $l_0\in\mathfrak{Primes}$ be a prime number which divides $N$.
Since $a\in (k^\times)_{\mathrm{div}}$, it holds that $\mu_{l_0^\infty}\subset k^\times$.
Therefore, $G_k^{l_0}$ is trivial.
\qed

\begin{Thm}\label{KF}
$k$ is Kummer-faithful if and only if $G_k\simeq\hat{\mathbb{Z}}$.
\end{Thm}

\prf

First, suppose that $k$ is Kummer-faithful.
By Proposition \ref{torally}, to show that $G_k\simeq\hat{\mathbb{Z}}$, it suffices to show that $G_k^p\simeq\mathbb{Z}_p$ (or, equivalently, $G_k^p$ is not trivial).
Suppose that $G_k^p$ is trivial.
Let $A$ be an abelian variety over $k$ whose $p$-rank is not zero.
By \cite[Proposition 2.2]{Ozeki}, there exists a finite extension $k'$ of $k$ such that $A(k')[p^\infty]=A(\overline{k})[p^\infty]$.
This contradicts the assumption that $k$ (and hence $k'$) is Kummer-faithful (cf. \cite[Proposition 2.4 (2)]{OT}).

Conversely, suppose that $G_k\simeq\hat{\mathbb{Z}}$.
Since $k$ is torally Kummer-faithful by Proposition \ref{torally}, it suffices to prove that $A(k)_{\mathrm{div}}=0$ for any abelian variety over $k$ (cf. \cite[Proposition 2.3]{OT}).
Let $k_0$ be a finite subfield of $k$ over which $A$ is defined.
Since $k$ is Galois over $k_0$ and $k_0$ is Kummer-faithful, by \cite[Proposition 2.4 (2)]{OT}, it suffices to show that $A(k)[l^\infty]$ is finite for any prime number $l$.
Since $G_k\simeq\hat{\mathbb{Z}}$, $k$ is a potential prime-to-$l$ extension of $k_0$ (for the definition of potential prime-to-$l$ extensions, see \cite[Definition 2.1]{Ozeki}).
Therefore, by \cite[Proposition 2.4]{Ozeki}, $A(k)[l^\infty]$ is finite, as desired.
\qed

\begin{Rem}
By Proposition \ref{torally} and Theorem \ref{KF}, algebraic extensions of finite fields of characteristic $p$ whose absolute Galois groups are isomorphic to $\hat{\mathbb{Z}}^{p'}$ are torally Kummer-faithful and not Kummer-faithful.
(For other examples of torally Kummer-faithful fields which are not Kummer-faithful, see \cite[Remark 1.5.3 (ii)]{Topics3}.)
\end{Rem}

\begin{Rem}
In \cite[\S 17.2]{Stix}, Stix gave examples of infinite algebraic extensions $\mathbb{F}$ of finite fields such that, for every abelian variety $A$ over $\mathbb{F}$, the Kummer map
\[A(\mathbb{F})\to H^1(G_{\mathbb{F}},\,T(A))\]
is injective, and that $G_{\mathbb{F}}\simeq\hat{\mathbb{Z}}$, where $G_{\mathbb{F}}$ is the absolute Galois group of $\mathbb{F}$.
\end{Rem}

\begin{Cor}\label{finGC}
Let $k_i$ be an algebraic extension of a finite field for $i=1,\,2$.
Suppose that the absolute Galois group of $k_i$ is isomorphic to $\hat{\mathbb{Z}}$ for $i=1,\,2$.
Let $X_i$ be a hyperbolic curve over $k_i$ and $\Pi_i$ either the \'{e}tale fundamental group $\pi_1(X_i)$ or the tame fundamental group $\pi_1^{\mathrm{tame}}(X_i)$ of $X_i$ for $i=1,\,2$ (for some choice of basepoint).
Write $\mathrm{Isom}(\Pi_1,\,\Pi_2)$ for the set of isomorphisms of profinite groups $\Pi_1\stackrel{\sim}{\to}\Pi_2$, $\mathrm{Isom}^{\mathrm{PT}}(\Pi_1,\,\Pi_2)\subset\mathrm{Isom}(\Pi_1,\,\Pi_2)$ (resp. $\mathrm{Isom}^{\mathrm{PT\text{-}GP}}(\Pi_1,\,\Pi_2)\subset\mathrm{Isom}(\Pi_1,\,\Pi_2)$) for the subset of point-theoretic (resp, point-theoretic and Galois-preserving) isomorphisms of profinite groups $\Pi_1\stackrel{\sim}{\to}\Pi_2$, $\mathrm{Isom}(X_1,\,X_2)$ for the set of isomorphisms of schemes $X_1\stackrel{\sim}{\to}X_2$, and $\mathrm{Inn}(\Pi_2)$ for the group of inner automorphisms of $\Pi_2$.
(For the definition of point-theoretic or Galois-preserving isomorphisms, see \cite[Definition 3.1]{Kummer}.)
Then the following assertions hold:

\begin{enumerate}
\item[(i)] Suppose that $\mathrm{Isom}(\Pi_1,\,\Pi_2)$ is not empty.
Then $\Pi_1=\pi_1(X_1)$ if and only if $\Pi_2=\pi_1(X_2)$.

\item[(ii)] Suppose that $\mathrm{Isom}(\Pi_1,\,\Pi_2)$ is not empty, and that $\Pi_i=\pi_1^{\mathrm{tame}}(X_i)$ for some (hence both of) $i=1,\,2$ (cf. (i)).
Then any element of $\mathrm{Isom}(\Pi_1,\,\Pi_2)$ is Galois-preserving (in particular, $\mathrm{Isom}^{\mathrm{PT}}(\Pi_1,\,\Pi_2)=\mathrm{Isom}^{\mathrm{PT\text{-}GP}}(\Pi_1,\,\Pi_2)$).

\item[(iii)] Suppose that either $X_1$ or $X_2$ is affine.
Then the natural map
\[\mathrm{Isom}(X_1,\,X_2)\to\mathrm{Isom}(\Pi_1,\,\Pi_2)/\mathrm{Inn}(\Pi_2)\]
determines a bijection
\[\mathrm{Isom}(X_1,\,X_2)\stackrel{\sim}{\to}\mathrm{Isom}^{\mathrm{PT\text{-}GP}}(\Pi_1,\,\Pi_2)/\mathrm{Inn}(\Pi_2).\]
\end{enumerate}
\end{Cor}

\prf

(i) follows from Proposition \ref{topfin} and Corollary \ref{type} (see also \cite[Theorem 3.4 (i)]{Kummer}).
(ii) follows from Proposition \ref{tame}.
(iii) follows from \cite[Theorem 3.4 (ii)]{Kummer}.
\qed

\

\

Let $d$ be a non-negative integer.
In the present paper, we shall say that a field $K$ is a {\it{$d$-dimensional local field}} if there is a chain of fields $K=K_d,\,K_{d-1},\,\cdots,\,K_1,\,K_0$ where $K_{i+1}$ is a complete discrete valuation field with residue field $K_i$ for $0\leq i\leq d-1$, and $K_0$ is of positive characteristic and algebraic over the prime field.
The field $K_{d-1}$ (resp. $K_0$) is said to be the {\it{first}} (resp. the {\it{last}}) {\it{residue field}} of $K$.
Moreover, we shall say that a field $K$ is a {\it{higher local field}} if $K$ is a $d$-dimensional local field for some non-negative integer $d$.
As in \cite[Remark 1.2]{HLF}, a chain of fields $K=K_d,\,K_{d-1},\,\cdots,\,K_1,\,K_0$ is unique up to isomorphism, and, in particular, the dimension $d$ of a higher local field $K$ is well-defined.

\begin{Cor}\label{GMLF}
Let $K$ be a higher local field of characteristic zero.
Suppose that the first residue field of $K$ is of positive characteristic and that the absolute Galois group of the last residue field of $K$ is isomorphic to $\hat{\mathbb{Z}}$.
Then $K$ is Kummer-faithful.
\end{Cor}

\prf

Immediate from Theorem \ref{KF} and \cite[Proposition 3.7]{HLF}.
\qed

\begin{Rem}
By Corollary \ref{GMLF}, mixed-characteristic complete discrete valuation fields whose residue fields are infinite algebraic extensions of finite fields with absolute Galois groups isomorphic to $\hat{\mathbb{Z}}$ are Kummer-faithful and generalized sub-$p$-adic for some $p\in\mathfrak{Primes}$, but not sub-$p$-adic.
Moreover, any higher local field satisfying the conditions in Corollary \ref{GMLF} whose last residue field is infinite is Kummer-faithful and not sub-$p$-adic (for any $p\in\mathfrak{Primes}$).
(For other examples of Kummer-faithful fields which are not sub-$p$-adic (for any $p\in\mathfrak{Primes}$), see \cite[Remark 1.5.4 (iii)]{Topics3} and \cite[\S 3]{OT}.)
\end{Rem}

\begin{Rem}
Let $K$ be a $d$-dimensional local field ($d\in\mathbb{Z}_{\geq 0}$).
Fix a chain of residue fields $K=K_d,\,K_{d-1},\,\cdots,\,K_1,\,K_0$, and let $G_{K_i}$ be the absolute Galois group of $K_i$ for $0\leq i\leq d$, and $I_{K_i}$ the kernel of the surjection $G_{K_i}\twoheadrightarrow G_{K_{i-1}}$ for $1\leq i\leq d$  (and set $G_K:=G_{K_d}$).
Similarly to \cite[\S 1]{FK}, one (and only one) of the following holds:
\begin{enumerate}
\item[(1)] $K$ is isomorphic to $K_0(\!(T_1)\!)\cdots(\!(T_d)\!)$.

\item[(2)] $K$ is isomorphic to $k(\!(T_1)\!)\cdots(\!(T_{d-1})\!)$, where $k$ is a field isomorphic to a finite extension of the quotient field of the Witt ring with coefficients in $K_0$.

\item[(3)] $K$ is isomorphic to a finite extension of $k\{\!\{T_1\}\!\}\cdots\{\!\{T_m\}\!\}(\!(T_{m+2})\!)\cdots(\!(T_d)\!)$, where $1\leq m\leq d-1$ and $k$ is a field isomorphic to a finite extension of the quotient field of the Witt ring with coefficients in $K_0$.
\end{enumerate}

In the following, consider the case where $d>0$ and $G_{K_0}\simeq\hat{\mathbb{Z}}$.
Set $p:=\mathrm{char}\,K_0$.
Note that we may characterize $p$ group-theoretically (in the sense of mono-anabelian reconstruction).
Indeed, if $G_K$ is topologically finitely generated, then $K_1$ is isomorphic to a $p$-adic local field (cf. \cite[Propositions 1.13, 2.1 and 2.2]{GMLF}), and therefore the group-theoreticity of $p$ follows from \cite[Theorem 2.15]{HLF}.
On the other hand, if $G_K$ is not topologically finitely generated, $p$ is characterized as the unique prime number $l$ such that $G_K^l$ is not topologically finitely generated.

Similarly to \cite[Proposition 2.1]{HLF}, we may determine group-theoretically (from $G_K$) which of the above three conditions $K$ satisfies:

\begin{enumerate}
\item[(i)] $K$ satisfies (1) if and only if $\mathrm{cd}_p\,G_K=1$.

\item[(ii)] $K$ satisfies (2) if and only if $\mathrm{cd}_p\,G_K\geq 2$, and, for any open normal subgroup $H$ of $G_K$, it holds that $H^1(H,\,\mathbb{Z}/p\mathbb{Z})$ is at most countable.

\item[(iii)] $K$ satisfies (3) if and only if $\mathrm{cd}_p\,G_K\geq 2$ and there exists an open normal subgroup $H$ of $G_K$ such that $H^1(H,\,\mathbb{Z}/p\mathbb{Z})$ is uncountable.
\end{enumerate}

Indeed, (i) is immediate from \cite[Proposition 2.2 (i)]{GMLF} and \cite[\S 0, Corollary]{K}.
Let $\mathcal{O}_K$ be the ring of integers of $K$ and $\mathfrak{M}_K$ the maximal ideal of $\mathcal{O}_K$.
Moreover, if $\mathrm{char}\,K=0$ and $\mathrm{char}\,K_{d-1}=p$, let $e_K$ be the absolute ramification index of $K$.

First, suppose that $K$ satisfies (2) and we shall prove that $H^1(H,\,\mathbb{Z}/p\mathbb{Z})$ is at most countable for any open normal subgroup $H$ of $G_K$.
By replacing $K$ by a finite extension, we assume that $G_K=H$.
We may also assume that $d=1$ since we have the following exact sequence for $2\leq i\leq d$:
\[\xymatrix{0 \ar[r] & H^1(G_{K_{i-1}},\,\mathbb{Z}/p\mathbb{Z}) \ar[r] & H^1(G_{K_i},\,\mathbb{Z}/p\mathbb{Z}) \ar[r] & H^1(I_{K_i},\,\mathbb{Z}/p\mathbb{Z})}\,(\simeq\mathbb{Z}/p\mathbb{Z}).\]
Moreover, by replacing $K$ by a finite extension if necessary, we may assume that $K$ contains a primitive $p$-th root of unity.
By Kummer theory, we have
\[H^1(G_K,\,\mathbb{Z}/p\mathbb{Z})\simeq\mathbb{Z}/p\mathbb{Z}\times (1+\mathfrak{M}_K)/(1+\mathfrak{M}_K)^p.\]
(Note that $K_0$ is perfect.)
Let $N$ be an integer satisfying $N>\dfrac{pe_K}{p-1}$.
Then we have $1+\mathfrak{M}_K^N\subset (1+\mathfrak{M}_K)^p$ and hence we obtain the following surjection:
\[(1+\mathfrak{M}_K)/(1+\mathfrak{M}_K^N)\twoheadrightarrow (1+\mathfrak{M}_K)/(1+\mathfrak{M}_K)^p.\]
Therefore, it suffices to show that $(1+\mathfrak{M}_K)/(1+\mathfrak{M}_K^N)$ is at most countable.
On the other hand, for any $j\in\mathbb{Z}_{>0}$, we have $(1+\mathfrak{M}_K^j)/(1+\mathfrak{M}_K^{j+1})\simeq K_0$.
Since $K_0$ is at most countable, $(1+\mathfrak{M}_K)/(1+\mathfrak{M}_K^N)$ is also at most countable, as desired.

Next, suppose that $K$ satisfies (3).
For simplicity, we assume that $d=2$ (we may treat the case where $d>2$ similarly).
By replacing $K$ by a finite extension if necessary, we assume that $K$ contains a primitive $p$-th root of unity.
By Kummer theory, we have
\[H^1(G_K,\,\mathbb{Z}/p\mathbb{Z})\simeq\mathbb{Z}/p\mathbb{Z}\times K_1^\times/(K_1^\times)^p\times (1+\mathfrak{M}_K)/(1+\mathfrak{M}_K)^p\twoheadrightarrow (1+\mathfrak{M}_K)/(1+\mathfrak{M}_K^2)\simeq K_1.\]
(Note that $(1+\mathfrak{M}_K)^p\subset 1+\mathfrak{M}_K^2$.)
However, $K_1\simeq K_0(\!(T)\!)$ is uncountable.
Therefore, $H^1(G_K,\,\mathbb{Z}/p\mathbb{Z})$ is also uncountable.

Moreover, the dimension $d$ of $K$ is also recovered group-theoretically (cf. \cite[\S 1]{EF} and \cite[Theorem 2.15 (i)]{HLF}):
\[d=\max_{\substack{l\in\mathfrak{Primes}\setminus\{p(G)\} \\ H\subset G}}\left(\dim_{\mathbb{Z}/l\mathbb{Z}}((H^l)^\mathrm{ab}\otimes_{\mathbb{Z}_l}\mathbb{Z}/l\mathbb{Z})\right)-1,\]
where $H$ runs through the set of open subgroups of $G_K$.

\end{Rem}

\begin{Cor}\label{GC}
Let $k_i$ be a higher local field of characteristic zero whose first residue field is of positive characteristic for $i=1,\,2$ (cf. the paragraph preceding Corollary \ref{GMLF}).
Suppose that the absolute Galois group of the last residue field of $k_i$ is isomorphic to $\hat{\mathbb{Z}}$ for $i=1,\,2$.
Let $X_i$ be a hyperbolic curve over $k_i$ and $\pi_1(X_i)$ the \'{e}tale fundamental group of $X_i$ for $i=1,\,2$ (for some choice of basepoint).
Write $\mathrm{Isom}(\pi_1(X_1),\,\pi_1(X_2))$ for the set of isomorphisms of profinite groups $\pi_1(X_1)\stackrel{\sim}{\to}\pi_1(X_2)$, $\mathrm{Isom}^{\mathrm{PT}}(\pi_1(X_1),\,\pi_1(X_2))\subset\mathrm{Isom}(\pi_1(X_1),\,\pi_1(X_2))$ for the subset of point-theoretic isomorphisms of profinite groups $\pi_1(X_1)\stackrel{\sim}{\to}\pi_1(X_2)$, $\mathrm{Isom}(X_1,\,X_2)$ for the set of isomorphisms of schemes $X_1\stackrel{\sim}{\to}X_2$, and $\mathrm{Inn}(\pi_1(X_2))$ for the group of inner automorphisms of $\pi_1(X_2)$.
(For the definition of point-theoretic isomorphisms, see \cite[Definition 3.1 (i)]{Kummer}.)
Suppose that either $X_1$ or $X_2$ is affine.
Then the natural map
\[\mathrm{Isom}(X_1,\,X_2)\to\mathrm{Isom}(\pi_1(X_1),\,\pi_1(X_2))/\mathrm{Inn}(\pi_1(X_2))\]
determines a bijection
\[\mathrm{Isom}(X_1,\,X_2)\stackrel{\sim}{\to}\mathrm{Isom}^{\mathrm{PT}}(\pi_1(X_1),\,\pi_1(X_2))/\mathrm{Inn}(\pi_1(X_2)).\]
\end{Cor}

\prf

By Corollary \ref{GMLF} and \cite[Theorem 3.4]{Kummer}, it suffices to prove that any isomorphism $\pi_1(X_1)\stackrel{\sim}{\to}\pi_1(X_2)$ of profinite groups is Galois-preserving.
This follows from \cite[Theorem D]{MT}.
\qed

\begin{Rem}
In the situation of Corollary \ref{GC}, if $k_i$ is a higher local field of dimension $1$ (of characteristic zero), $k_i$ is generalized sub-$p$-adic for some $p\in\mathfrak{Primes}$ (for $i=1,\,2$).
Therefore, the ``relative Isom-version'' of the Grothendieck conjecture for hyperbolic curves over these fields is affirmative (cf. \cite[Theorem 4.12]{sur}).
However, Corollary \ref{GC} gives a ``semi-absolute'' result, which is unknown even for hyperbolic curves over $p$-adic local fields in general (without point-theoretic assumption).
\end{Rem}

\begin{Rem}
A ``weak version'' of the Grothendieck conjecture for hyperbolic curves of genus $0$ over subfields of finitely generated extensions of mixed-characteristic higher local fields is already known (cf. \cite[Theorem A]{MT}).
\end{Rem}

\

\

\end{document}